\documentclass{amsart}

\usepackage[T1]{fontenc}
\usepackage[utf8]{inputenc}
\usepackage[USenglish]{babel}
\usepackage{textcase}

\usepackage[
	bookmarks=true,
	plainpages=false,
	linktocpage,
	colorlinks=true,
	citecolor=green!80!black,
	linkcolor=red!70!black,
	filecolor=magenta,
	urlcolor=magenta,
	breaklinks,
	pdfauthor={Martin Winter},
]{hyperref}

\usepackage{amsmath,amsthm,amssymb}
\usepackage{calc, mathtools} 

\usepackage{colortbl,color,xcolor} 
\usepackage{centernot} 
\usepackage{array} 
\usepackage{enumitem,moreenum} 
\usepackage{cite} 
\usepackage[nameinlink,capitalize,noabbrev]{cleveref}
\usepackage{nicefrac}

\usepackage{tikz-cd} 

\usepackage[font=small,labelfont=bf]{caption}

\newcommand{\ZZ}{\mathbb{Z}}    
\newcommand{\RR}{\mathbb{R}}    
                   
\newcommand{\Tsymb}{\top}
\newcommand{\T}{^{\Tsymb}}

\def\^#1{^{(#1)}}
\def\s^#1{^{\smash{(#1)}}}

\newcommand*{\horzbar}{\rule[.5ex]{2.5ex}{0.4pt}}

\newcommand{\wideword}[1]{\quad\text{#1}\quad}
\newcommand{\wideand}{\wideword{and}}

\def\:{\colon}

\newcommand{\mylabel}{\upshape(\textit{\roman*})}
\newenvironment{myenumerate}{\begin{enumerate}[label=\mylabel]}{\end{enumerate}}

\newcommand{\freespace}{\kern.07em} 
\newcommand{\free}{\freespace\cdot\freespace} 

\newcommand{\defeq}{\mathrel{\mathop:}=}                         

\newcommand{\enquote}[1]{``#1''}                                 

\theoremstyle{plain}  
\newtheorem{theorem}{Theorem}[section]
\newtheorem{corollary}[theorem]{Corollary}
\newtheorem{lemma}[theorem]{Lemma}
\newtheorem{proposition}[theorem]{Proposition}

\theoremstyle{definition} 
\newtheorem{definition}[theorem]{Definition}

\newtheorem{construction}[theorem]{Construction}

\crefname{theorem}{Theorem}{Theorems}
\crefname{proposition}{Proposition}{Propositions}
\crefname{lemma}{Lemma}{Lemmas}
\crefname{corollary}{Corollary}{Corollaries}
\crefname{remark}{Remark}{Remarks}
\crefname{example}{Example}{Examples}
\crefname{definition}{Definition}{Definitions}
\crefname{problem}{Problem}{Problems}
\crefname{observation}{Observation}{Observation}
\crefname{construction}{Construction}{Construction}

\DeclareMathOperator{\Aut}{Aut}

\DeclareMathOperator{\Ortho}{O}
\DeclareMathOperator{\SOrtho}{SO}
\DeclareMathOperator{\GL}{GL}

\DeclareMathOperator{\Span}{span}

\DeclareMathOperator{\rank}{rank}

\DeclareMathOperator{\Id}{Id}
\DeclareMathOperator{\id}{id}

\DeclareMathOperator{\Fix}{Fix}

\DeclareMathOperator{\Sym}{Sym}
\DeclareMathOperator{\Cov}{Cov} 		
\DeclareMathOperator{\Perm}{Perm}   	

\let\<=\langle
\let\>=\rangle
\let\x=\times

\def\...{...}
\newcommand{\shortStyle}{\textit}
\newcommand{\ie}{\shortStyle{i.e.,}}

\newcommand{\eg}{\shortStyle{e.g.}}
\newcommand{\Eg}{\shortStyle{E.g.}}
\newcommand{\wrt}{\shortStyle{w.r.t.}}

\newcommand{\resp}{resp.}

\makeatletter
\renewcommand*{\eqref}[1]{%
  \hyperref[{#1}]{\textup{\tagform@{\ref*{#1}}}}%
}
\makeatother

\numberwithin{equation}{section}

\def\img{img/}

\begin{document}


\expandafter\title
{Geometry and topology of symmetric point arrangements}
		
\author[M. Winter]{Martin Winter}
\address{Faculty of Mathematics, University of Technology, 09107 Chemnitz, Germany}
\email{martin.winter@mathematik.tu-chemnitz.de\newline\rule{0pt}{1.5cm}\includegraphics[scale=0.7]{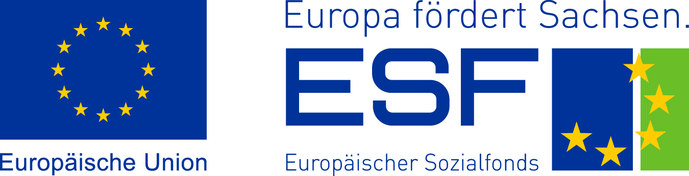}}
	
\subjclass[2010]{20B25, 52C25}
\keywords{symmetries, point arrangements, arrangement space, permutation groups}
		
\date{\today}
\begin{abstract}
We investigate point arrangements $v_i\in\RR^d,i\in \{1,\...,n \}$ with certain prescribed symmetries.
The \emph{arrangement space} of $v$ is the column span of the matrix in which the $v_i$ are the rows. 
We characterize properties of $v$ in terms of the arrangement space, \eg\ we characterize whether an arrangement possesses certain symmetries or whether it can be continuously deformed into another arrangement while preserving symmetry in the process.
We show that whether a symmetric arrangement can be continuously deformed into its mirror image depends non-trivially on several factors, \eg\ the decomposition of its representation into irreducible constituents, and whether we are in even or odd dimensions.


\end{abstract}

\maketitle


\section{Introduction}
\label{sec:introduction}

A \emph{point arrange\-ment} (or just arrangement) is a finite family of points \mbox{$v_i\in\RR^d$}, $i\in N :=\{1,\...,n\}$. 
An arrangement is symmetric \wrt\ some permutation group $\Gamma\subseteq\Sym(N)$, if any permutation $\phi\in\Gamma$ can be realized on the points via some isometry of the ambient space (formally defined in \cref{sec:symmetric_arrangements}). 
Central to our treatment of~point arrangements is the notion of \emph{arrangement spaces}.
To introduce these, define the matrix
$$M :=
\begin{pmatrix}
	\;\horzbar\!\!\!\! & v_1\T & \!\!\!\!\horzbar\;\; \\
	& \vdots & \\
	\;\horzbar\!\!\!\! & v_n\T & \!\!\!\!\horzbar\;\;
\end{pmatrix}\in\RR^{n\times d}$$
in which the $v_i$ are the rows. 
We call $M$ the \emph{arrangement matrix} (or just matrix) of $v$.
The arrangement space $U:=\Span M\subseteq\RR^n$ is then defined as the column span of $M$.
Arrangements with the same arrangement space are called \emph{equivalent}.

The definition of arrangement space is motivated by a recurring idea in geometry, and despite its simplicity has some interesting and non-trivial applications. 
To our knowledge, no common name (or no name at all) was given to this concept so far. 
We list some of its applications.


The notion of equivalence (\ie\ having the same arrangement space) has a direct geometric interpretation: two arrangements of the same dimension are equivalent, if and only if they are related by an invertible linear transformation (this~is~a~well-known fact from linear algebra: two matrices $M,\bar M\in\RR^{n\x d}$ have the same column span if and only if $M=\bar M T$ for some invertible $T\in\GL(\RR^d)$; see also \cref{res:spherical_point_arrangement_subspace}).
A prominent use of arrangement spaces is therefore in areas where one cares less about the exact positioning of the points, and more about linear, affine or convex dependencies between them, as \eg\ in the study of point configurations and oriented matroids (see \eg\ \cite[Section 6.3]{ziegler2012lectures}).
In this setting, a subspace $U\subseteq\RR^n$ can be thought of as defining an equivalence class of point arrangements \wrt\ invertible linear transformations.
A representative of this class is obtained as follows:
choose vectors $u_1,\...,u_d\in U$ that span $U$, and define the matrix $M:=(u_1,\...,u_d)\in\RR^{n\x d}$ in which the $u_i$ are the columns.
The rows of $M$ are then a point arrangement with matrix $M$ and arrangement space $U$.
%
%

Properties of an arrangement that are invariant under invertible linear transformations are already determined by the arrangement space. 
For example, the rank of an arrangement $\rank v:=\rank\{v_1,\...,v_n\}$ equals the dimension of the arrangement space $U$. In particular, there is always a full-dimensional arrangement $v$ (\ie\ $\rank v=d$) for a given arrangement space. 
Likewise, \mbox{whether an~arrangement~is~a} linear transformation of a rational arrangement $v_i\in\mathbb{Q}^d$ or a 01-arrangement $v_i\in$ $\{0,1\}^d$ is determined by the arrangement space.

For studying point configurations and polytopes in higher dimensions, there exists the notion of \emph{Gale duality} (see \cite[Section 6]{ziegler2012lectures}).
The construction of the Gale dual is usually considered as quite artificial and technical.
The following diagram gives a natural construction via arrangement spaces:
%
\begin{center}
\begin{tikzcd}
v   \arrow[r, "\text{Gale}", dashed, mapsto]
    \arrow[d, mapsto]& \bar v 
    \arrow[d, mapsfrom] \\U\arrow[r, mapsto]& U^{\mathrlap\bot}
\end{tikzcd}
\end{center}
%
In other words: $v$ and $\bar v$ are Gale duals if and only if their arrangement spaces are orthogonal complements of each other.
Many properties are shared between $U$ and its orthogonal complement (\eg, possessing a rational basis, being invariant \wrt\ some orthogonal
 representation, etc.), and so it is clear that many properties are expressed by $v$ and its Gale dual likewise.

We close with an application to \emph{linear matroids} (see, \eg\ \cite{gordon2012matroids} for a general reference, and Section 6 in particular). 
The linear matroid of an arrangement $v$ is given by its independent sets, that is, subsets $I\subseteq N$ for which $\{v_i\in\RR^n \mid i\in I\}$ is linearly independent. 
The same matroid can be obtained from the arrangement space $U$ of $v$.
This time, as independent sets we choose the subsets $I\subseteq N$ for which $\Span\{e_i\in\RR^n\mid i\in I\}$ has trivial intersection with $U^\bot$ ($e_i\in\RR^n$ is the $i$-th standard basis vector). 
Call this the matroid of $U$. 
There exists the notion of a dual matroid, and it is not hard to see that the dual of the matroid of $U$ is the matroid of $U^\bot$. 
As noted previously, this is exactly the linear matroid of the Gale dual of $v$. 
This gives a short proof of the fact that the dual of a linear matroid is linear.

The general objective of this paper is to apply the idea of arrangement spaces to symmetric point arrangements.
While the \enquote{symmetry} in \enquote{symmetric arrangement} is easily defined with the language of representation theory, we demonstrate that the study of their topology involves intricacies that are not obvious from just studying group representations.

This paper is the first in a row of three papers (see the outlook in \cref{sec:future}),~all of which deal with symmetries of discrete geometric objects: points in Euclidean space (this paper), graph embeddings \cite{winter2020symmetric} and convex polytopes \cite{winter2020eigenpolytopes}.
Our terminology is therefore motivated from a geometric point of view, 
but is still strongly related to well-established terminology in other areas, most notably, (finite) frame theory (see \eg\cite{waldron2018introduction,casazza2013introduction}).
For example, in frame theory, point arrangements are known as \emph{finite frames}, and their associated arrangement matrices are known as \emph{analysis operators}.
We later introduce spherical, normalized and symmetric arrangements which in frame theory are referred to as \emph{tight frames}, \emph{Parseval frames} and \emph{group frames} \cite{waldron2018group} respectively.
Our results and the general approach of this paper originated within the mental image of points in Euclidean space, and we encourage the reader to think of point arrangements as such clouds of points.
Nevertheless, one can equivalently view all of this from the perspective of frame theory.
We support this second perspective by providing, for certain basic theorems in \cref{sec:point_arrangements,sec:symmetric_arrangements}, references to equivalent versions in the literature of frame theory. 
Still, we always included a proof in our terminology in order to be self-contained.

Finally, also in the geometric context our treatment of point arrangements with symmetries has to be distinguished from a list of similarly flavored ideas, most notably, \emph{orbit polytopes} \cite{babai1977symmetry, friese2016affine} and \emph{bar-joint-frameworks} with symmetry constraints \cite{schulze2010symmetry,malestein2014generic,nixon2016symmetry}.
Our notions are more general in that we do not require all points to lie on the boundary of their convex hull, nor do we discuss distances constraints between the points.

\subsection*{Overview}

\Cref{sec:preliminaries} contains a short proof of a version of Schur's lemma~rele\-vant for this paper.

From \cref{sec:point_arrangements} on our investigations focus on so-called \emph{normalized} arrangements. We show that the arrangement space determines these up to orthogonal transformations, and thus, determines their metric properties. 

\Cref{sec:symmetric_arrangements} gives a formal definition of \enquote{symmetric arrangement}, \ie\ \emph{$\Gamma$-arrange\-ments} for some permutation group $\Gamma\subseteq\Sym(N)$. 
We show that the arrangement spaces of symmetric arrangements are invariant subspaces of $\RR^n$, and that this property characterizes symmetric arrangements.

\Cref{sec:rigid_flexible} takes up a major part of this paper and discusses the topological aspects of point arrangements.
Here we apply the developed techniques to answer questions of the following kind: given two $\Gamma$-arrangements $v,\bar v$, is it possible to continuously deform $v$ into $\bar v$ while preserving the symmetry in the process? 
As will turn out, the arrangement space plays a role in answering these questions.
For example, we shall obtain the following result: if $v$ and $\bar v$ are non-equivalent irreducible arrangements of odd dimension, and their arrangement spaces are non-orthogonal, then $v$ can be continuously deformed into $\bar v$ in the sense explained above.

\section{Preliminary facts}
\label{sec:preliminaries}

We formulate and prove a version of Schur's lemma for orthogonal representations. 
Recall that for two $\Gamma$-representations $T,\bar T\:\Gamma\to\Ortho(\RR^d)$, a map $R\in \RR^{d\x d}$ with $R\, T_\phi =\bar T_\phi R$ for all $\phi\in\Gamma$ is called \emph{equivariant}.

\begin{theorem}[Schur's lemma]\label{res:schur}
Given two linear representations $T,\bar T\:\Gamma\to\Ortho(\RR^d)$. 
If at least one of these is irreducible, then every equivariant map between them is of the form $\alpha R$ for some orthogonal $ R\in\Ortho(\RR^d)$.
\end{theorem}

It is convenient to briefly repeat the proof for this specific form.
For re-usability, we first prove the following lemma:

\begin{lemma}\label{res:schur_prep}
If $T\:\Gamma\to\Ortho(\RR^d)$ is an irreducible representation and $R\in\RR^{d\x d}$ is symmetric and commutes with $T_\phi$ for all $\phi\in\Gamma$, then $R=\alpha\Id$ for some $\alpha\in\RR$.
\begin{proof}
Since $R$ commutes with $T_\phi$ for all $\phi\in\Gamma$, each eigenspace of $R$ is preserved by each $T_\phi$. So, each such preserved eigenspace is an invariant subspace of $T$. But since $T$ is irreducible, $T$ has no non-trivial invariant subspace. Hence, $R$ must have a single eigenspace to eigenvalue, say, $\alpha\in\RR$. And since $R$ is symmetric, it is diagonalizable, and we have $R=\alpha \Id$.
\end{proof}
\end{lemma}

We proceed with the proof of Schur's lemma:

\begin{proof}[Proof of Schur's lemma (\cref{res:schur})]
Without loss of generality, assume that $T$ is irreducible.
Let then $R$ be an equivariant map between $T$ and $\bar T$, that is, $R\, T_\phi=\bar T_\phi R$, or $R=\bar T_\phi R\, T_\phi\T\!$ for all $\phi\in\Gamma$. We compute $R\T\! R$:
$$R\T\! R=T_\phi R \T\!\underbrace{\bar T_\phi\T\!\bar T_\phi}_{=\,\Id} R\, T_\phi\T\! = T_\phi R\T\!R\, T_\phi\T\!.$$
Rearranging shows that the symmetric matrix $R\T\! R$ commutes with $T_\phi$ for all $\phi\in\Gamma$, and by \cref{res:schur_prep} we can conclude $R\T\! R=\tilde \alpha \Id$ for some $\tilde\alpha\in\RR$. Since $R\T\! R$ is positive semi-definite, we have $\tilde\alpha\ge 0$, and we can choose $\alpha\in\{\pm\tilde\alpha^{1/2}\}$ so that $\alpha^{-1}R$ is orthogonal.
\end{proof}

Recall that two representations are said to be \emph{isomorphic}, if there is an \emph{invertible} equivariant map between them. We immediately obtain the following corollary:

\begin{corollary}
Given two linear representations $T,\bar T\:\Gamma\to\Ortho(\RR^d)$. 
If at least one of them is irreducible and there exists a non-zero equivariant map between them, then they are isomorphic.
\end{corollary}

Isomorphic representations are considered identical for all practical purposes, as they are related by a simple change of basis. 
In particular, one is irreducible if and only if the other one is.

\section{Normalized arrangements}
\label{sec:point_arrangements}

Our investigations are motivated from geometry, in which metric properties of point arrangements (such as angles and lengths) are the focus of interest.
However, our goal is also to study point arrangements mainly via their arrangement spaces, and these determine the arrangement only up to invertible linear transformations.
Such general transformations do not respect the metric properties of an arrangement (\eg\ a square and a rhombus have the same arrangement space, but are metrically quite different).

We therefore decided on a compromise: while not strictly necessary for the purpose of this paper, we will often focus our investigation on a sub-class of arrangements for which any two equivalent arrangements are the same up to an \emph{orthogonal} transformation.
We shall see later that this comes with no loss of generality, but at times, is more convenient to work with.
These sub-classes are as follows:

\begin{definition}
Let $v$ be an arrangement with matrix $M$. 
\begin{myenumerate}
	\item $v$ is called \emph{spherical}, if $M\T\!M=\alpha \Id$ for some $\alpha\in\RR\setminus\{0\}$. 
	\item $v$ is called \emph{normalized}, if $M\T\!M=\Id$.
\end{myenumerate}
\end{definition}

The term \enquote{spherical} stems from the interpretation of $M\T\!M$ as (a scaled version of) the \emph{covariance matrix}:
%
\begin{equation}\label{eq:covariance}
\Cov(v):=\frac1{n-1}\sum_{i=1}^n v_i v_i\T\!\quad\implies\quad M\T\!M=(n-1)\Cov(v).
\end{equation}
Spherical (and normalized) arrangements are always full-dimensional:
$$ \rank v=\rank M = \rank (M\T\!M)=\rank (\alpha \Id) = d.$$
Finally, as advertised before, normalized arrangements are especially desirable for the following property (a version for frames can be found in \cite[Proposition 3.2]{waldron2018introduction}):

\begin{theorem}\label{res:spherical_point_arrangement_subspace}
Two normalized arrangements $v,\bar v$ are equivalent if and only if they are related by an orthogonal transformation $T\in\Ortho(\RR^d)$, in particular, both are $d$-dimensional. 
The orthogonal transformation is given by
%
$$ T=\bar M\T\! M,$$
where $M,\bar M$ are the matrices of $v$ and $\bar v$ respectively.
\begin{proof}
Let $U,\bar U\subseteq\RR^n$ be the arrangement spaces of $v$ and $\bar v$ respectively.

If $v=T\bar v$ for some invertible linear transformation $T\in\GL(\RR^d)$, then this can be written as $M=\bar M T\T\!$. This proves the first direction via
$$U=\Span M=\Span(\bar MT\T)=\Span \bar M = \bar U.$$
%

For the other direction, assume that $v$ and $\bar v$ are equivalent, that is, $U=\bar U$.
Since normalized arrangements are full-dimensional, $\rank v = \dim U = \dim \bar U = \rank \bar v$ alrady proves that both are of the same dimension. 
We have to check that $T:=\bar M\T\! M$ satisfies $MT\T\! = \bar M$:
$$MT\T\! = M M\T\! \bar M \overset{(*)}= \bar M,$$
where in $(*)$ we used that $MM\T\!$ acts as ortho-projector onto $U$. Since the columns of $\bar M$ are in $\bar U = U$, $MM\T\!$ acts as identity.
It remains to show that $T$ is orthogonal:
$$T T\T\!=\bar M\T\! M M\T\! \bar M = \bar M\T\! \bar M=\Id,$$
where $MM\T\!$ was dropped for the same reason as above.
%
\end{proof}
\end{theorem}


From now on we often focus on normalized (or spherical) arrangements (and we state explicitly when we do so).
This choice emphasizes our intent and comes with no loss of generality:
the columns of the arrangement matrix of a normalized arrangement form an orthonormal system (ONS) in $\RR^n$. 
For a general full-dimensional arrangement $v$, by simply ortho-normalizing the columns of its arrangement matrix, we obtain a normalized arrangement that is equivalent to $v$
(in the language of frames, the existence of this so-called \emph{canonical tight frame} is shown in \cite[Theorem 3.4]{waldron2018introduction}).

%

Our intent is further reflected in the definitions and results of \cref{sec:symmetric_arrangements,sec:rigid_flexible}:
in \cref{sec:symmetric_arrangements} we consider only symmetries by orthogonal transformations (rather than invertible linear transformations), and in \cref{sec:rigid_flexible} we consider only deformations that keep the arrangements normalized (rather than full-dimensional).
Again, neither comes with a loss of generality (see especially \cref{res:irreducible_is_spherical}).




\section{Symmetric arrangements}
\label{sec:symmetric_arrangements}

We now consider point arrangements with certain geometric symmetries. The symmetries will be prescribed by some permutation group $\Gamma\subseteq\Sym(N)$.
To get an idea for what this means, consider the following: if there were a permutation $\phi\in\Gamma$ that exchanges $1$ and $2$, we then would require the existence of an orthogonal transformation in $\RR^d$ that exchanges points $v_1$ and $v_2$ (it can happen that $v_1=v_2$, and then any transformation achieves this; but that is the boring case).

We will investigate how point arrangements with such symmetries can be constructed and classified.
In fact, we will see that symmetric arrangements can be fully characterized by their arrangement space.


\begin{definition}\label{def:Gamma_arrangement}
A point arrangement $v$ is called \emph{$\Gamma$-arrangement}, if there exists a representation $T\:\Gamma\to\Ortho(\RR^d)$ with
\begin{equation}\label{eq:symmetry_condition}
T_\phi v_i = v_{\phi(i)},\quad\text{for all $i\in N$ and $\phi\in\Gamma$.}
\end{equation}
$T$ is called a \emph{representation} of $v$.
\end{definition}

For us, all representations will be defined over $\RR$. 
In particular, a representation is considered as \emph{irreducible}, if it is \emph{real} irreducible.
Also, as mentioned previously, we focus on \emph{orthogonal} representations because our motivation stems from a metric context. This is no loss of generality.

Condition \eqref{eq:symmetry_condition} in the above definition can be stated using the matrix $M$ of $v$: 
\begin{equation}\label{eq:PM_MT}
P_\phi M=MT_\phi,\quad\text{for all $\phi\in \Gamma$},
\end{equation}
where $P_\phi\in\Perm(n)$ is the permutation matrix that permutes the components of a vector in $\RR^n$ according to $\phi$.
This statement is equivalent to $P_\phi\T\! M= M T_\phi\T\!$ (which can be easily verified in this form by writing out the matrix $M$), but
we will mostly apply it in the form \eqref{eq:PM_MT}.
%
%
We also see that $M$ and $M\T\!$ are now equivariant maps between $T$ and the permutation matrix representation $\phi\mapsto P_\phi$ in $\RR^n$.

The representation of a $\Gamma$-arrangement does not have to be uniquely determined by its set of points. 
If, however, $v$ is full-dimensional, the representation is indeed unique: 
if $\Span\{v_1,\...,v_n\}=\RR^d$, we can find a basis of $\RR^d$ of points, say, $v_{i_1},\...,v_{i_d}\in\RR^d$. 
For each $\phi\in\Gamma$, the transformation $T_\phi\in\Ortho(\RR^d)$ is then already uniquely determined by
$$T_\phi v_{i_j}=v_{\phi(i_j)},\quad\text{for all $j\in\{1,\...,d\}$}.$$
We will therefore often speak of \emph{the} representation of $v$.


\begin{lemma}\label{res:equivariant_map_MM}
If $v,\bar v$ are two $\Gamma$-arrangements with matrices $M, \bar M$, then $\bar M\T\! M$ is an equivariant map between the representations of $v$ and $\bar v$.
\begin{proof}
This follows from \eqref{eq:PM_MT} and its equivalent form: let $T,\bar T$ be representations of $v$ and $\bar v$, then
$$\bar T_\phi  \bar M\T\! M = \bar M\T\! P_\phi M = \bar M\T\! M T_\phi.$$
\end{proof}
\end{lemma}


At the moment, we still discuss general symmetric arrangements, that is, not necessarily spherical or normalized.
However, recall that if normalized arrangements $v$ and $\bar v$ are related by an orthogonal transformation, then \cref{res:spherical_point_arrangement_subspace} states that this transformation is $\bar M\T\! M$. 
This transformation is then invertible (because orthogonal). 
So, if $v,\bar v$ are $\Gamma$-arrangements, the same transformation $\bar M\T\! M$ is also an equivariant map between the representations $T$ and $\bar T$, and we can conclude that the representations are isomorphic.

\begin{corollary}\label{res:equivalent_impies_isomophic}
Equivalent normalized $\Gamma$-arrangements have isomorphic representations.
\end{corollary}

%
%

The converse of \cref{res:equivalent_impies_isomophic} is \emph{not} true: if two (normalized) arrangements have isomorphic representations, they might still be non-equivalent. We will discuss this further in \cref{sec:rigid_flexible}. For now, note the freedom of choice for $v_1$ in the following construction, which might result in non-equivalent arrangements:

\begin{construction}
We construct a $\Gamma$-arrangement with prescribed representation that acts transitively on its point set.

Choose a transitive permutation group $\Gamma\subseteq\Sym(N)$ and a representation $T\:\Gamma\to\Ortho(\RR^d)$ of it.
Let $\Gamma_1\subset\Gamma$ be the stabilizer of $1\in N$.
Consider the following subspace:
$$\Fix(T,\Gamma_1):= \{x\in\RR^d\mid T_\phi x= x,\text{ for all $\phi\in\Gamma_1$}\}\subseteq\RR^d.$$
Choose a $v_1\in \Fix(T,\Gamma_1)$. Since $\Gamma$ is transitive, for any $i\in N$ there is a $\phi_i\in\Gamma$ with $\phi_i(1)=i$.
We then define the desired arrangement as follows:
$$v_i:= T_{\phi_i} v_1,\quad\text{for all $i\in N$}.$$
This will turn out to be a $\Gamma$-arrangement with representation $T$. 


If $\phi\in\Gamma$ maps $i\in N$ to $j$, then $\psi:= \phi_j^{-1}\circ \phi\circ \phi_i\in\Gamma_1$, and
$$T_\phi v_i =T_{\phi_j} T_{\psi} v_1 = T_{\phi_j} v_1 = v_j,$$
since $v_1\in\Fix(T,\Gamma_1)$ and thus is fixed by $T_{\psi}$.
This shows that, in particular, the construction is independent of the choice of the $\phi_i\in\Gamma$.
We then have that $T$ is a $\Gamma$-representation as needed in \cref{def:Gamma_arrangement}, and thus that $v$ is a $\Gamma$-arrangement with representation $T$.
\end{construction}

This construction essentially proves
that $\Fix(T, \Gamma_1)$ is (canonically) isomorphic to the set of equivariant maps between $T$ and the permutation representation of $\Gamma$.

Note further that in the above construction, a necessary condition for the existence of non-zero~arrangements with representation $T$ is $\Fix(T,\Gamma_1)\not=\{0\}$.
A necessary condition for the existence of multiple non-equivalent $\Gamma$-arrangements with the prescribed representation is $\dim \Fix(T,\Gamma_1)\ge 2$. 
However, this is not sufficient.

Symmetric arrangements inherit a notion of irreducibility from their representations:
a $\Gamma$-arrangement shall be called \emph{irreducible} if its representation is irreducible (over $\RR$).
It shall be called \emph{reducible} otherwise. 
There are equivalent ways to characterize reducibility in geometric terms (see \cref{res:redicible_formuations} and \cref{fig:arrangement}). 
To describe these, we require the following construction:

\begin{figure}
\centering
\includegraphics[width=0.7\textwidth]{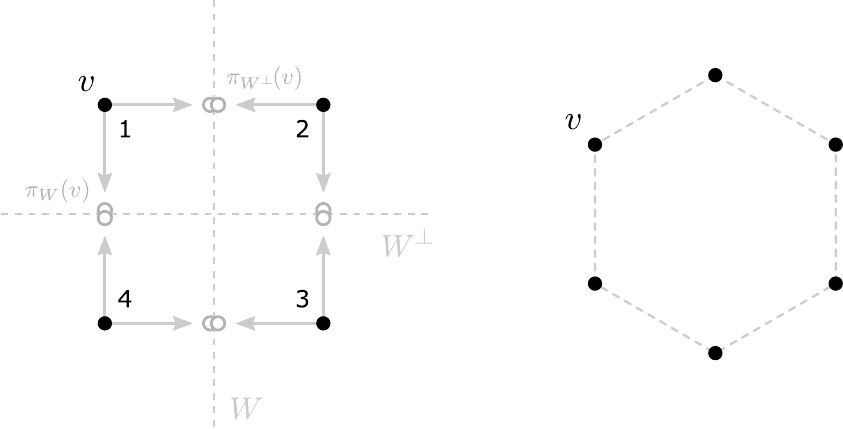}
\caption{Left: a $\Gamma$-arrangement with group $\Gamma\simeq\ZZ_2\x\ZZ_2$ acting on four points via the permutations $\{\id, (12)(34),(13)(24),(14)(23)\}$. The arrangement is reducible, as it can be projected onto $W$ and $W^\bot$ without loosing symmetry (note that $90^\circ$ rotations are not a prescribed symmetry). Right: an irreducible $\Gamma$-arrangement with $\Gamma\simeq D_6$ (symmetry group of the regular hexagon) on six points. This arrangement cannot be projected to a lower dimension while keeping rotational symmetry.}
\label{fig:arrangement}
\end{figure}

\begin{definition}
Given arrangements $v$ and $\bar v$ of dimension $d$ and $\bar d$ respectively, the \emph{direct sum arrangement} $v\x\bar v$ is defined by
$$(v\oplus \bar v)_i:=(v_i,\bar v_i)\in\RR^{d+\bar d},\quad\text{for all $i\in N$}.$$
\end{definition}

If $M,\bar M$ are the matrices of the arrangements $v$ and $\bar v$ respectively, the direct sum arrangement has the matrix $(M,\bar M)$, \ie\ the columns of $M$ and $\bar M$ are joined to a single matrix. 
If $U,\bar U\subseteq\RR^n$ are the arrangement spaces of $v$ and $\bar v$, then the direct sum has the arrangement space $U+ \bar U$.

Direct sum arrangements can be used to obtain new $\Gamma$-arrangements from old ones:
it is not hard too see that sums of $\Gamma$-arrangements are again $\Gamma$-arrangements. 
If $T,\bar T$ are representations of $v$ and $\bar v$, then $v\oplus \bar v$ has representation
\begin{equation}\label{eq:product_representation}
\phi\quad\mapsto\quad (T\oplus\bar T)_\phi :=\begin{pmatrix}
\,T_\phi & \\ & \!\!\bar T_\phi
\end{pmatrix}.
\end{equation}

We list some geometric characterizations of being reducible:

\begin{lemma}\label{res:redicible_formuations}
Let $v$ be a non-zero $\Gamma$-arrangement. The following are equivalent:
\begin{myenumerate}
\item \label{it:reducible} 
The arrangement $v$ is reducible.
\item \label{it:product} 
The arrangement $v$ is (equivalent to) the direct sum of two $\Gamma$-arrangements of smaller dimension.
\item \label{it:projection} 
There exists a proper subspace $W\subset\RR^n$\!, so that the projections $\pi_W(v)$ and $\pi_{W^\bot\!}(v)$ are $\Gamma$-arrangements with the same representation as $v$.
\end{myenumerate}
\begin{proof}
Assume \ref{it:reducible}, that is, $v$ is reducible. 
By this, $T$ is reducible, and there exists a proper non-zero $T$-invariant subspace $W\subset\RR^d$. 
Consider the projection $\bar v:=\pi_W(v)$ onto $W$. The ortho-projector $\pi_W$ commutes with $T_\phi$ for all $\phi\in\Gamma$, hence $\bar v$ is a $\Gamma$-arrangement with representation $T$:
$$T_\phi \bar v_i = T_\phi \pi_W(v_i) = \pi_W(T_\phi v_i)=\pi_W (v_{\phi(i)}) = \bar v_{\phi(i)}.$$
The orthogonal complement $W^\bot$ is $T$-invariant as well and the same reasoning applies to the projection $\pi_{W^\bot\!}(v)$. 
This proves \ref{it:projection}.

We now assume \ref{it:projection}. 
By equivalence, we can assume 
$$W=\Span\{e_1,\...,e_{\delta}\}\wideand W^\bot=\Span\{e_{\delta+1},\...,e_d\}$$ 
(where $e_i$ denote the $i$-th standard basis vectors). 
If $M=(u_1,\...,u_d)$ is the matrix of $v$, the projected arrangements have matrices
$$M_W = (u_1,\..., u_\delta)\wideand M_{W^\bot}=(u_{\delta+1},\...,u_d).$$
We can consider these projections as arrangements in $\RR^\delta$ \resp\ $\RR^{d-\delta}$. The join of the matrices is obviously $M$, and thus the direct sum of the projections is $v$.
This proves \ref{it:product}.


Finally, assume \ref{it:product}, that is, $v = v\^1\oplus v\^2$ is a $d$-dimensional $\Gamma$-arrangement with representation $T$, and decomposes into $d_i$-dimensional $\Gamma$-arrangements $v\s^i$ with representation $T\s^i$. 
By \eqref{eq:product_representation} we have $T=T\s^1\oplus T\s^2$.
%
%
Then $W:=\Span\{e_1,\...,e_{d_1}\}$ is $T$-invariant, and hence $T$ is reducible and \ref{it:reducible} holds.

\end{proof}
\end{lemma}

The irreducible $\Gamma$-arrangements are the elementary building blocks of general $\Gamma$-arrangements: each $\Gamma$-arrangement $v$ can be written as a direct sum
$$v= v\^1\oplus \cdots\oplus  v\^K$$
with irreducible $\Gamma$-arrangements $v\s^k$, for $k\in\{1,\...,K\}$. 
This decomposition into irreducible constituents might not be unique. 

As elaborated in \cref{sec:point_arrangements}, we are especially interested in spherical and normalized arrangements. 
General $\Gamma$-arrange\-ments can be far from spherical, but there is no concern in the case of irreducible $\Gamma$-arrangements (see also \cite[Theorem 10.5]{waldron2018introduction} for a version in terms of groups frames):

\begin{proposition}\label{res:irreducible_is_spherical}
Irreducible $\Gamma$-arrangements are spherical.
\begin{proof}
Let $v$ be an irreducible $\Gamma$-arrangement with representation $T$.
By \cref{res:equivariant_map_MM} (with setting $v=\bar v$), we see that $M\T\! M$ is an equivariant map between $T$ and $T$, \ie\ it commutes with $T$.
Then, $M\T\! M=\alpha\Id$ follows from \cref{res:schur_prep}.
\end{proof}
\end{proposition}

In the reducible case, consider two normalized $\Gamma$-arrangements $v$ and $\bar v$. Then
$$\alpha v\oplus \beta \bar v$$
is a $\Gamma$-arrangement, but to be spherical we necessarily need $|\alpha|=|\beta|$ (this is not sufficient, though).

The following two theorems give a characterization of $\Gamma$-arrangements in terms of their arrangement spaces.
For this, note that $\Gamma$ acts on $\RR^n$ via $\phi\mapsto P_\phi$, hence we can speak of $\Gamma$-invariant subspaces of $\RR^n$.
With the following results we answer the question of classification of symmetric arrangements (by shifting the problem~to a classification of $\Gamma$-invariant subspaces).

\begin{theorem}\label{res:Gamma_arrangement_has_gamma_invariant_subspace}
The arrangement space of a $\Gamma$-arrangement is $\Gamma$-invariant.
\begin{proof}
Let $v$ be a $\Gamma$-arrangement with matrix $M$, arrangement space $U\subseteq\RR^n$ and representation $T$.
For every $\phi\in\Gamma$ holds
\vspace{-0.3em}
$$P_\phi U = \Span (P_\phi M) = \Span(MT_\phi)=\Span M = U.$$
%
Hence $U$ is $\Gamma$-invariant.
\end{proof}
\end{theorem} 

The converse of \cref{res:Gamma_arrangement_has_gamma_invariant_subspace} is true for spherical arrangements:

\begin{theorem}\label{res:Gamma_invariant_subspace_implies_gamma_arrangement}
A spherical arrangement with $\Gamma$-invariant arrangement space is a $\Gamma$-arrange\-ment.
Moreover, its representation can be chosen as
\begin{equation}\label{eq:T_MPM}
T_\phi:= M\T\!P_\phi M,\quad\text{for all $\phi\in\Gamma$},
\end{equation}
where $M$ is the matrix of the arrangement.
\begin{proof}
Let $v$ be a spherical arrangement, $M$ its matrix, and $U\subseteq\RR^n$ its $\Gamma$-invariant arrangement space.
The arrangement $v$ is a $\Gamma$-arrangement if and only $\alpha v$ is~a~$\Gamma$-arrangement for every $\alpha\in\RR\setminus\{0\}$, and both have the same arrangement space.
Hence, we can assume that $v$ is normalized.

We show that $v$ is a $\Gamma$-arrangement by proving that \eqref{eq:T_MPM} is indeed the desired representation.
Let $\phi,\psi\in\Gamma$ and observe
$$T_\phi T_\psi = (M\T\! P_\phi M)(M\T\! P_\psi M)=(M\T\! P_\phi)(MM\T\!)(P_\psi M).$$
Since $v$ is normalized, $MM\T\!$ acts as ortho-projector onto $U$. 
Since $\Span M=U$ is invariant \wrt\ $P_\psi$, the columns of $P_\psi M$ are in $U$ again. 
Thus, $MM\T\!$ acts as identity to its right and can be dropped. 
We obtain 
$$T_\phi T_\psi = M\T\! P_\phi P_\psi M = M\T\! P_{\phi\circ \psi} M = T_{\phi\circ\psi}.$$
This shows that $T$ is a linear representation.
Equivalently, one shows that $T_\phi\T=T_{\phi^{\smash{-1}}}$, \ie\ that $T_\phi$ is orthogonal.

It remains to prove $T_\phi M\T=M\T\! P_\phi$ for all $\phi\in\Gamma$:
$$T_\phi M\T = (M\T\! P_\phi M) M\T = \big(MM\T( P_{\phi^{\smash{-1}}} M)\big)\T.$$
Again, $MM\T\!$ acts to its right (on the columns of $P_{\phi^{\smash{-1}}}M$) as ortho-projector, hence as identity, and can be dropped. We obtain $T_\phi M\T = (P_{\phi^{\smash{-1}}} M)\T = M\T\! P_\phi$ and we are done.
\end{proof}
\end{theorem} 

The preceding theorem is not necessarily true for non-spherical arrangements. \Eg\ the vertices of a rectangle \resp\ square share the same arrangement space (they are linear transformations of each other), but only one is a $\Gamma$-arrangement for $\Gamma\simeq D_4$ (the symmetry group of the square).

\Cref{res:Gamma_arrangement_has_gamma_invariant_subspace} and \cref{res:Gamma_invariant_subspace_implies_gamma_arrangement} together make clear that studying spherical $\Gamma$-arrangements is essentially just studying $\Gamma$-invariant subspaces of $\RR^n$. 

\begin{corollary}\label{res:charcterization_Gamma_arrangement}
A spherical arrangement is a $\Gamma$-arrangement if and only if its arrangement space is a $\Gamma$-invariant subspace of $\RR^n$.
\end{corollary}

It is in general a non-trivial task to determine the $\Gamma$-invariant subspaces of $\RR^n$, and hence it is non-trivial to determine the $\Gamma$-arrangements. 
But as soon as a $\Gamma$-invariant subspace is known, we can obtain a $\Gamma$-arrangement as a representative of the equivalence class described by this subspace.

We close this section by proving that a full-dimensional $\Gamma$-arrangement is irreducible if and only if its arrangement space is irreducible as a $\Gamma$-invariant subspace of $\RR^n$.
The same is not true for $\Gamma$-arrangements that are not of full dimension. 
Such arrangements are always reducible, regardless of their arrangement space, since $\Span v = \Span\{v_1,\...,v_n\}$ is an invariant subspace of $T$.

\begin{proposition}\label{res:irreducible_arrangement_and_subspace}
A full-dimensional $\Gamma$-arrangement $v$ is irreducible if and only if its arrangement space $U\subseteq\RR^n$ is irreducible as a $\Gamma$-invariant subspace of $\RR^n$.
\begin{proof}
Let $M$ be the matrix of $v$, and $T$ its representation. 
Since $v$ is full-dimen\-sio\-nal, we can interpret $M$ as an invertible linear map $M\:\RR^d \to U$.

By \cref{res:Gamma_arrangement_has_gamma_invariant_subspace}, the map $\phi\mapsto P_\phi|_U$ is a $\Gamma$-representation on $U$ (since $U$ is $\Gamma$-invariant).
By \eqref{eq:PM_MT} $M$ is then an invertible equivariant map between $T$ and $\phi\mapsto P_\phi|_U$.
Hence, these representations are isomorphic, and, in particular, one is irreducible if and only if the other one is.

The arrangement $v$ is irreducible if and only if $T$ is. The subspace $U$ is irreducible if and only of $\phi\mapsto P_\phi|_U$ is. These statements about irreducibility are hence linked since $T$ and $\phi\mapsto P_\phi|_U$ are isomorphic.
\end{proof}
\end{proposition}


\section{Topology of symmetric arrangements}
\label{sec:rigid_flexible}

From this section on, for the rest of the paper, we study the topology of symmetric point arrangements. 
Observe that the set of all $d$-dimensional arrangements indexed by $N$ carries the structure of an $nd$-dimensional vector space, and hence is equipped with a natural topology.

\begin{figure}
\centering
\includegraphics[width=0.8\textwidth]{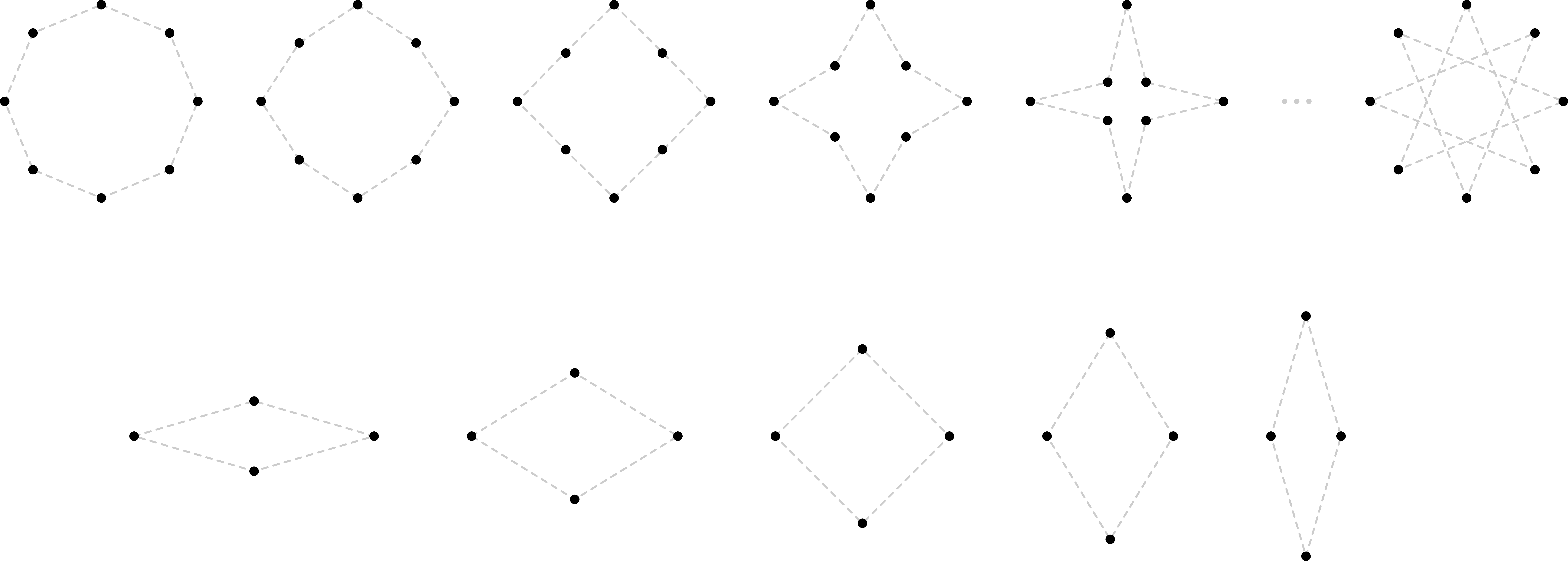}
\caption{An intuitive visualization of the concepts of deformation and flexibility. The shown transitions preserve symmetries defined by groups $\Gamma\cong D_4$ on eight points, and $\Gamma\cong \ZZ_2\x\ZZ_2$ on four points respectively (neither of these groups acts transitively). While the upper one stays spherical (and is a deformation as defined in the text), the bottom one does not. In fact, restricting to spherical arrangements during a deformation makes the bottom arrangement rigid (only the arrangement as a square is possible).}
\label{fig:flexible}
\end{figure}

We shall be concerned with the following type of questions:
suppose we have two (normalized) $\Gamma$-arrangements $v$ and $\bar v$. 
We try to continuously deform $v$ into $\bar v$.
This is certainly possible if we are allowed to move the points freely.
We therefore restrict to deformations, in which every intermediate step is itself a $\Gamma$-arrangement, \ie\ the symmetry is preserved.
There are still trivial solutions, as every arrangement can be continuously deformed to the zero-arrangement (which is always symmetric) via $[0,1]\ni t\mapsto tv$.
For that reason, we additionally impose the constraint that the arrangement stays normalized during the deformation process.
See \cref{fig:flexible} for a visualization.

We can achieve \enquote{deformations} in trivial cases, \eg\ when $\bar v = T v$ for some orientation preserving orthogonal transformation $T\in\SOrtho(\RR^d)$.
We use the fact that $\SOrtho(\RR^d)$ is path-connected in $\RR^{d\x d}$ to construct a continuous transition between $v$ and $Tv$.
Of course, this is not really a deformation in the intuitive sense, but more a continuous re-orientation.
The situation is less obvious if we consider orientation reversing transformations, \ie\ $\bar v=Tv$ with $T\in\Ortho(\RR^d)\setminus\SOrtho(\RR^d)$.
In fact, the answer for when an arrangement can be continuously deformed into its mirror image depends on the dimension of the arrangement.

We finally consider \enquote{proper deformations} between non-equivalent $\Gamma$-arrange\-ments, \ie\ where the arrangements are not just re-orientations of each other.
In particular, we ask for which $\Gamma$-arrangements such a non-trivial deformation is possible at all.
This gives a notion of \emph{rigidity} for $\Gamma$-arrangements.
We characterize rigidity of an arrangement in terms of its arrangement space.


As motivated above, we consider deformations that preserve symmetry and normalization. 
For that matter, for each $\Gamma\subseteq\Sym(N)$, we consider the set
$$\mathcal S_d(\Gamma) := \{\,\text{$d$-dimensional normalized $\Gamma$-arrangements on $N$}\,\},$$
equipped with the subspace topology.
%

\begin{definition}
A $d$-dimensional \emph{$\Gamma$-deformation} (or just \emph{deformation}) is a curve (that is, a continuous function) $v(\free)\:[0,1]\to\mathcal S_d(\Gamma)$. 

$\Gamma$-arrangements $v,\bar v\in\mathcal S_d(\Gamma)$ are \emph{deformation equivalent}, if there exists a $\Gamma$-deformation $v(\free)$ with $v(0)=v$ and $v(1)=\bar v$.
We then also say more intuitively that $v$ can be \emph{deformed} into $\bar v$.
\end{definition}

\begin{figure}
\centering
\includegraphics[width=0.6\textwidth]{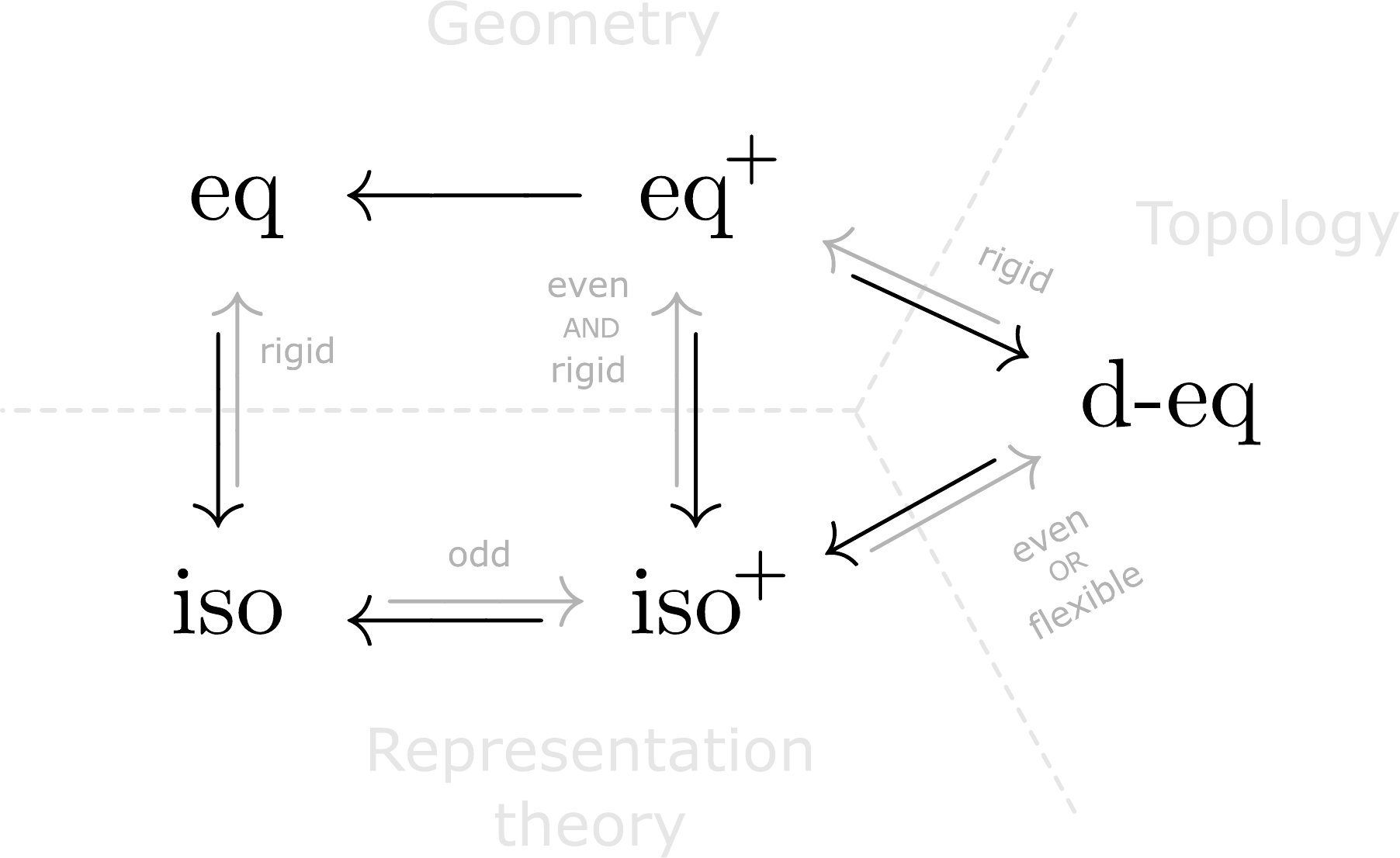}
\caption{Implications between the five relations defined on $\mathcal S_d(\Gamma)$. The black arrows indicate unconditional implications. The gray arrows indicate conditional implications, and the conditions are written next to the arrow. The conditions are necessary. The conditions \enquote{rigid/flexible} mean that at least one (and then both) of the involved arrangements is rigid/flexible. The conditions \enquote{odd/even} mean that the arrangements live in odd/even dimensional space.}
\label{fig:relations}
\end{figure}

We already reasoned that for our question it will be important to distinguish between orientation preserving and orientation reversing transformations. 
We therefore define the following more specific versions of equivalence and isomorphy of representations: 
\begin{itemize}
\item
Two arrangements are \emph{positively equivalent} if they are related by a linear transformation of positive determinant.
\item
Two representations are \emph{positively isomorphic} if there exists an equivariant map between these with positive determinant.
\end{itemize}

In total, we now defined five relations on $\mathcal S_d(\Gamma)$, each of which is easily seen to be an equivalence relation: equivalence (eq), positive equivalence ($\mathrm{eq}^+$), representation-iso\-mor\-phy (iso; \ie\ having isomorphic representations), positive representation-iso\-mor\-phy ($\mathrm{iso}^+$; \ie\ having positively isomorphic representations) and deformation equivalence (d-eq).

A major goal of this section is to prove the \emph{eleven} implications presented in \cref{fig:relations}.\footnote{Keep this figure at hand, as we refer back to it quite frequently. You may mark the already proven arrows to keep track of the progress.}
Some of the implications are conditional. 
For these, we also prove that the conditions are necessary.
We will focus here on the irreducible arrangements in $\mathcal S_d(\Gamma)$.
Most results can be modified to work for reducible arrangements as well, but we will not pursue this. 

The following implications in \cref{fig:relations} are trivially seen to be true, or where already discussed in the introduction to this section:
$$
\mathrm{eq}^+\! \longrightarrow \mathrm{eq},\qquad
\mathrm{iso}^+\! \longrightarrow \mathrm{iso},\qquad
\mathrm{eq}^+\! \longrightarrow \text{d-eq}.
$$
Further, if normalized $\Gamma$-arrangements $v,\bar v$ are equivalent, then the relating transformation (which is $\bar M\T\! M$) is simultaneously an invertible equivariant map between their representations (see also \cref{res:equivalent_impies_isomophic}). This proves
$$
\mathrm{eq} \longrightarrow \mathrm{iso},\qquad
\mathrm{eq}^+\! \longrightarrow \mathrm{iso}^+\!.
$$

We next aim to prove a necessary condition for two arrangements to be deformation equivalent. The corresponding statement can be found in \cref{fig:relations} in the form of the unconditional implication $\text{d-eq}\longrightarrow\mathrm{iso}^+$.  In other words: 
to be deformation equivalent, it is necessary to have positively isomorphic representations.
This seems to be quite natural: the set of (pair-wise non-isomorphic) representations of $\Gamma$ is finite. The discrete nature of this set does not work well with the continuity of deformations.
We make this precise: for a $\Gamma$-representation $T\:\Gamma\to\Ortho(\RR^d)$ let
$$\mathcal S_d^+(T):=\{v\in\mathcal S_d(\Gamma)\mid \text{$v$ has a representation positively isomorphic to $T$}\}.$$
The sets $\mathcal S_d^+(T)$ are the equivalence classes in $\mathcal S_d(\Gamma)$ under the relation of positive representation-isomorphy.
We prove the following topological result:

\begin{theorem}\label{res:S_d_plus_open}
$\mathcal S_d^+(T)$ is an open subset of $\mathcal S_d(\Gamma)$.
\end{theorem}

For what follows, it is convenient to introduce the following function:
$$\det(v,\bar v):=\det(\bar M\T\! M),$$
where $v,\bar v$ are arrangements, and $M,\bar M$ their matrices.
We list some of its properties:


%

\begin{lemma}[Properties of $\det(\free,\free)$]\label{res:det_properties}
The following holds:
\begin{myenumerate}
\item \label{it:continuity}
$\det(\free,\free)$ is continuous in both arguments.
\item \label{it:normalized_1}
If normalized arrangements $v,\bar v$ are equivalent, then they have $\det(v,\bar v)=\pm1$, and if they are positively equivalent, then they have $\det(v,\bar v)=1$.
\item \label{it:orthogonal_0}
If two arrangements $v$ and $\bar v$ have orthogonal arrangement spaces, then they also have $\det(v,\bar v)=0$.
\item \label{it:0_orthogonal}
If $v,\bar v$ are irreducible $\Gamma$-arrangements with $\det(v,\bar v)=0$, then their arrangement spaces are orthogonal.
\item \label{it:isomorphy_no0}
If $\Gamma$-arrangements $v,\bar v$ have $\det(v,\bar v)\not=0$, then they have isomorphic representations, and if $\det(v,\bar v)>0$, then they have positively isomorphic representations.
\end{myenumerate}
\begin{proof}
Statement \ref{it:continuity} is directly seen from the definition.

In the following, let $v,\bar v$ be arrangements with arrangements matrices $M$ and $\bar M$ respectively.

If $v$ and $\bar v$ are normalized and equivalent, then $\bar M\T\! M$ is the orthogonal transformation between these. This already shows $$\det(v,\bar v)=\det(\bar M\T\! M)=\pm 1.$$
If the arrangements are positively isomorphic, then the transformation $\bar M\T\! M$ has positive determinant. This proves \ref{it:normalized_1}.

If the arrangement spaces $\Span M$ and $\Span \bar M$ are orthogonal, then each column of $M$ is orthogonal on each column of $\bar M$. This then shows $\bar M\T\! M=0$, in particular $\det(\bar M\T\! M)=0$, and hence \ref{it:orthogonal_0}.

For the rest of the proof, assume that $v$ and $\bar v$ are $\Gamma$-arrangements. Recall that then $\bar M\T \! M$ is an equivariant map between any representations of $v$ and $\bar v$.

To prove \ref{it:0_orthogonal}, assume that $v,\bar v$ are irreducible. By Schur's lemma $ \bar M\T\! M=\alpha R$ for some $R\in\Ortho(\RR^d)$. 
Hence, if $\det(v,\bar v)=\det(\bar M\T\! M)=0$, then only because $\alpha=0\implies\bar M\T\! M=0$, \ie\ $\Span \bar M\perp\Span M$.

On the other hand, to show \ref{it:isomorphy_no0}, observe that $\det(\bar M\T\! M)\not=0$ implies that $\bar M\T\! M$ is invertible. Since there is then an invertible equivariant map, the representations are isomorphic. 
If $\det(\bar M\T\! M)>0$, they are positively isomorphic.
\end{proof}
\end{lemma}

We have the following immediate consequence of the above properties, by simply combining \ref{it:0_orthogonal} with \ref{it:isomorphy_no0}:

\begin{corollary}\label{res:non_ortho_implies_isomorphy}
If irreducible arrangements $v,\bar v\in\mathcal S_d(\Gamma)$ have non-orthogonal arrangement spaces, then they have isomorphic representations.
\end{corollary}

We proceed with the proof of the theorem.

\begin{proof}[Proof of \cref{res:S_d_plus_open}]
We show that every $v\in \mathcal S_d^+(T)$ has an open neighborhood in $\mathcal S_d(\Gamma)$ which is completely contained in $\mathcal S_d^{\smash+}(T)$. 

Fix an arrangement $v\in \mathcal S_d^+(T)$, and consider the set
$$S:=\{\bar v\in\mathcal S_d(\Gamma) \mid \det(v,\bar v)>0\}.$$
$S$ is the pre-image of the open set $\RR_{>0}$ under a continuous function $\det(v,\free)$, hence open.
Also, $v\in S$, and so $S$ is an open neighborhood of $v$ in $\mathcal S_d(\Gamma$).

It remains to show that $S\subseteq\mathcal S_d^+(T)$. 
Choose a $\bar v\in S$, in particular, $\det(v,\bar v)>0$.
By \cref{res:det_properties}, the representations of $v$ and $\bar v$ are positively isomorphic, and $\bar v\in\mathcal S_d^+(T)$.
\end{proof}

We now see that the equivalence classes of $\mathcal S_d(\Gamma)$ under $\mathrm{iso}^+\!$ provide a decom\-position of $\mathcal S_d(\Gamma)$ into disjoint open sets.
It is well-known that points (in our case, arrangements) in distinct open components cannot be connected by a continuous path.
This shows $\text{d-eq}\longrightarrow\mathrm{iso}^+$.

\begin{corollary}\label{res:d_eq_implies_iso_plus}
Deformation equivalent arrangements $v,\bar v\in\mathcal S_d(\Gamma)$ have positively isomorphic representations.
\end{corollary}

This now motivates the question whether the $\mathcal S_d^+(T)$ are already the \emph{path-connected} components of $\mathcal S_d(\Gamma)$, \ie\ whether positive isomorphy of the representations already suffices to imply deformation equivalence. 

Surprisingly, the only case that we cannot immediately deal with, is, when the arrangements are equivalent.
We start by constructing a deformation between non-equivalent arrangements.
From now on, we restrict to irreducible arrangements.

\begin{proposition}\label{res:non_equivalent_deformation}
If two non-equivalent irreducible arrangements $v,\bar v\in\mathcal S_d(\Gamma)$ have positively isomorphic representations, then they are deformation equivalent.
\begin{proof}
The idea is as follows: first re-orient $\bar v$ so that the representations of $v$ and $\bar v$ become identical.
Second, make a linear transition between the arrangements.
Finally, make sure that all intermediate arrangements are normalized $\Gamma$-arrangements.

If $T,\bar T\:\Gamma\to\Ortho(\RR^d)$ are the positively isomorphic representations of $v$ and $\bar v$~res\-pectively, and since $v,\bar v$ are irreducible, Schur's lemma (\cref{res:schur}) ensures that there is an orthogonal equivariant map
 $R\in\SOrtho(\RR^d)$ of positive determinant between these, \ie\ $T_\phi=R\,\bar T_\phi R\T\!$ for all $\phi\in\Gamma$.
Then, $\bar v$ and $\tilde v:= R\,\bar v$ are positively equivalent, hence deformation equivalent by the already shown implication $\mathrm{eq}^+\!\longrightarrow \text{d-eq}$.
The arrangement $\tilde v$ is a $\Gamma$-arrangement with representation $T$:
$$T_\phi \tilde v_i = R\,\bar T_\phi R\T\! R\, \bar v_i = R\, \bar T_\phi \bar v_i = R\,\bar v_{\phi(i)} = \tilde v_{\phi(i)}.$$
It remains to construct a deformation from $v$ to $\tilde v$. Concatenating this with the deformation between $\tilde v$ and $\bar v$ then proves the statement. 

We define the continuous map (not necessarily a deformation, because not necessarily normalized)
$$t\quad\mapsto\quad v'(t):= (1-t) v + t \tilde v.$$
Remark that $v'(t)\not=0$ for all $t\in[0,1]$. 
This is clear for $t\in\{0,1\}$. If there is a $t\in(0,1)$ with $v'(t)=0$, then we can rearrange this to $v = t/(1-t)\cdot \tilde v$.
But then, $v$ is just a scaled version of $\tilde v$, in particular, equivalent to $\tilde v$, and by this also equivalent to $\bar v$.
This is in contradiction to the assumption that $v$ and $\bar v$ are non-equivalent.

Since both $v$ and $\tilde v$ are $\Gamma$-arrangements with representation $T$, so is $v'(t)$ for all $t\in[0,1]$.
In particular, $v'(t)$ is irreducible, and by \cref{res:irreducible_is_spherical} $v'(t)$ is spherical.
This means that if $M'(t)$ is the matrix of $v'(t)$, then $M'(t)\T\! M'(t)=\alpha(t) \Id$ for some continuous function $\alpha\:[0,1]\to\RR\setminus\{0\}$.
We finally define
$$t\quad \mapsto \quad v(t):= \alpha(t)^{-1/2} \cdot v'(t),$$
which is continuous.
Then all $v(t)$ are normalized, and $v(\free)$ is the desired deformation between $v$ and $\tilde v$.
\end{proof}
\end{proposition}

The remaining case of equivalent arrangements with positively isomorphic representations turns out to be surprisingly non-trivial.
In fact, we will observe that the answer depends on the parity of the dimension: positive representation-isomorphy suffices in even dimensions to imply the existence of a deformation, but not always in odd dimension.
The following lemma explains where the differences arise.

\begin{lemma}\label{res:representations_parity}
Let $T,\bar T\:\Gamma\to\Ortho(\RR^d)$ be two isomorphic irreducible representations. 
Depending on the parity of $d$, the following holds:
\begin{myenumerate}
\item \label{it:odd_case}
if $d$ is odd, then $T$ and $\bar T$ are already positively isomorphic.
\item  \label{it:even_case}
if $d$ is even, then either all equivariant maps between $T$ and $\bar T$ have non-positive determinant, or all of them have non-negative determinant.
\end{myenumerate}
\begin{proof}
Let $R, \bar R$ be two invertible equivariant maps between $T$ and $\bar T$. 

Assume that $d$ is odd. 
As a multiple of the identity, $-\Id$ commutes with $T_\phi$ for all $\phi\in\Gamma$.
Thus, also $-R$ is equivariant:
$$-R\, T_\phi = -\Id R\, T_\phi = -\Id \bar T_\phi R = \bar T_\phi (-\Id) R=\bar T_\phi (-R).$$
Because $d$ is odd, we have $\det(-R)=-\det(R)$, and either $R$ or $-R$ is an equivariant map of positive determinant, which proves \ref{it:odd_case}.

Now assume that $d$ is even.
We observe that $R\T\!\bar R$ commutes with $T_\phi$ for all permutations $\phi\in\Gamma$:
$$R\T\!\bar R\, T_\phi = R\T\! \bar T_\phi \bar R = T_\phi R\T\! \bar R.$$
By \cref{res:schur_prep}, $R\T\!\bar R=\alpha \Id$ for some $\alpha\in\RR$. 
But since $d=2\delta$ is even, and $R$ and $\bar R$ are invertible, we have that 
$$\det(R)\det(\bar R)=\det(R\T\!\bar R)=\det(\alpha \Id)=\alpha^{2\delta}$$
is positive, no matter the sign of $\alpha$.
Thus, $\det (R)$ and $\det (\bar R)$ have the same sign, which proves \ref{it:even_case}.
\end{proof}
\end{lemma}

In particular, the above lemma proves the implication $\mathrm{iso}\longrightarrow\mathrm {iso}^+$ under the assumption of odd dimension, and also shows that there are counterexamples in even dimensions.

We now show that in even dimensions and for irreducible $T$, the $\mathcal S_d^+(T)$ are indeed the path-connected components of $\mathcal S_d(\Gamma)$.

\begin{theorem}\label{res:even_positively_isomorphic_iff_deformation_equivalent}
Let $d$ be even. Irreducible arrangements $v,\bar v\in\mathcal S_d(\Gamma)$ are deformation equivalent if and only if their representations are positively isomorphic.
\begin{proof}
We have already shown one direction (see \cref{res:d_eq_implies_iso_plus}).
It remains to show that $\mathcal S_d(T)$ is path-connected.

If $v,\bar v\in\mathcal S_d(T)$ are non-equivalent, then a deformation exists by \cref{res:non_equivalent_deformation}.
We therefore can assume that they are equivalent. 
Then they are related by the orthogonal transformation $\bar M\T\! M\in\Ortho(\RR^d)$, where $M,\bar M$ are the matrices of $v$ and $\bar v$.
We know that $\bar M\T\! M$ is also an equivariant map between the representations of $v$ and $\bar v$.
Since the representations are positively isomorphic, and $d$ is even, by \cref{res:representations_parity}  all equivariant maps between the representations must have positive determinant.
In particular, $\det(\bar M\T\! M)>0$, and $v$ and $\bar v$ are positively equivalent.
We then apply the already proven implication $\mathrm{eq}^+\longrightarrow \text{d-eq}$.
\end{proof}
\end{theorem}

This proves one part of the implication $\mathrm{iso}^+\longrightarrow \text{d-eq}$ in \cref{fig:relations} (the one under the assumption of even dimension). The other part will be proven when we discuss odd dimensions.

For the remainder of this section, let $\tau\in\Ortho(\RR^d)\setminus\SOrtho(\RR^d)$ be some orientation reversing transformation, \ie\ $\det(\tau)=-1$. We call $\tau$ a reflection or mirror operation.
In even dimensions, reflections cannot be reversed by a continuous deformation:

\begin{corollary}
If $d$ is even and $v\in\mathcal S_d(\Gamma)$ an irreducible arrangement, then $v$ and $\tau v$ are not deformation equivalent.
\begin{proof}
Since $\tau$ is the transformation between $v$ and $\tau v$, it is also an equivariant map between their representations.
Since $d$ is even, by \cref{res:representations_parity} all equivariant maps between the representations must have a negative determinant.
Thus, the representations are not positively isomorphic, and $v$ and $\tau v$ are not deformation equivalent by \cref{res:even_positively_isomorphic_iff_deformation_equivalent}.
\end{proof}
\end{corollary}

In odd dimensions, the situation is more diverse. 
In particular, whether $\mathcal S_d^+(T)$ is path-connected depends on the group $\Gamma$ and the particular representation $T$.
The conditions can be stated in terms of rigidity:

\begin{definition}
A $\Gamma$-arrangement $v\in\mathcal S_d(\Gamma)$ is \emph{rigid} if it cannot be deformed into any non-equivalent arrangement. Otherwise it is called \emph{flexible}.
\end{definition}

In other words, every deformation between rigid arrangement is actually a continuous re-orientation.
And since a continuous re-orientation can only deform an arrangement into a positively equivalent arrangement, this proves the implication $\text{d-eq}\longrightarrow \mathrm{eq}^+$ (under the assumption of rigidity, see \cref{fig:relations}).

There is a straightforward way to characterize rigidity in terms of arrangement spaces as follows: an arrangement $v$ is rigid, if any deformation starting in $v$ has constant arrangement space during the deformation. This simply follows from \cref{res:spherical_point_arrangement_subspace}.

We list further characterizations of rigidity:

\begin{theorem}\label{res:rigidity_characterizations}
Let $v\in\mathcal S_d(\Gamma)$ be irreducible with arrangement space $U\subseteq\RR^n$ and representation $T$. The following are equivalent:
\begin{myenumerate}
\item \label{it:flexible} The arrangement $v$ is flexible.
\item \label{it:representation} There exists a non-equivalent arrangement $\bar v\in\mathcal S_d(\Gamma)$ with a representation isomorphic to $T$.
\item \label{it:subspace} There exists a $\Gamma$-invariant subspace $\bar U\subseteq \RR^n$ different from $U$ that is non-orthogonal to $U$.
\end{myenumerate}
\begin{proof}
The main tools of the proof will be the properties of $\det(\free,\free)$ (see \cref{res:det_properties}). 
We will prove $$\text{\ref{it:flexible} $\Longrightarrow$ \ref{it:subspace} $\Longrightarrow$ \ref{it:representation} $\Longrightarrow$ \ref{it:flexible}}.$$

Assume \ref{it:flexible} that is, $v$ is flexible, and there is a deformation $v(\free)$ with $v(0)=v$, and $v(1)$ some non-equivalent arrangement. 
%
If $\det(v,v(1))\not= 0$, then \ref{it:subspace} already follows from \cref{res:det_properties}. 
We therefore assume $\det(v,v(1))=0$. 
Now, consider the function $$\gamma(t):=\det(v,v(t)),\quad\text{for $t\in[0,1]$},$$ which is continuous in $t$, and satisfies $\gamma(0)=1$ and $\gamma(1)=0$. 
By the intermediate value theorem, there is a $t\in(0,1)$ with, say, $\gamma(t)=\det(v,v(t))=1/2$.
Since then $\det(v,v(t))\not=0$, we have that the arrangement spaces of $v$ and $v(t)$ are non-orthogonal by \cref{res:det_properties}, and because of $\det(v,v(t))<1$ it follows from the same result that the arrangement spaces are also distinct. 
Lastly, the arrangement space of $v(t)$ is $\Gamma$-invariant by \cref{res:Gamma_arrangement_has_gamma_invariant_subspace}, and hence we proved \ref{it:subspace}.

Now, assume \ref{it:subspace}, and choose a representative $\bar v\in\mathcal S_d(\Gamma)$ with arrangement space $\bar U$. 
Since $U$ and $\bar U$ are non-orthogonal, by \cref{res:non_ortho_implies_isomorphy} $v$ and $\bar v$ have isomorphic representations.
Since their arrangement spaces are also distinct, they are not equivalent, and we obtain \ref{it:representation}.

Finally, assume \ref{it:representation}, and let $\bar v\in\mathcal S_d(\Gamma)$ be the non-equivalent arrangement with isomorphic representation. 
In particular, $\bar v$ is irreducible.
Either $\bar v$ or $\tau\bar v$ (which is also non-equivalent to $v$) has a representation that is \emph{positively} isomorphic to $T$, hence we can assume that this holds for $\bar v$.
Then, there is a deformation between $v$ and $\bar v$ by \cref{res:non_equivalent_deformation}.
Since $v$ and $\bar v$ are non-equivalent, we see that $v$ must be flexible, which gives \ref{it:flexible}.
\end{proof}
\end{theorem}

From this result we learn that in order for a certain permutation group $\Gamma\subseteq\Sym(N)$ to allow for flexible $\Gamma$-arrangements, it is necessary that the $\Gamma$-represen\-ta\-tion $\phi\mapsto P_\phi$ has an irreducible constituent of multiplicity at least two. 
Such a constituent gives rise to several non-equivalent $\Gamma$-arrangements with isomorphic representations, hence flexibility by \cref{res:rigidity_characterizations}. Conversely, if all irreducible constituents are distinct, then all $\Gamma$-arrangements are rigid.

We now go back to the question of when $\mathcal S_d^+(T)$ is path-connected for $d$ odd.
While in even dimensions $v$ and $\tau v$ are not deformation-equivalent because they cannot have positively isomorphic representations, the same reasoning does not hold up in odd dimensions.
Here, isomorphy already implies positive isomorphy by \cref{res:representations_parity}.

We will now see that surprisingly, an odd-dimensional irreducible arrangement can be deformed into its mirror image if and only if it can be deformed at all, \ie\ if it is not rigid.

\begin{proposition}\label{res:odd_mirror_deformation}
If $d$ is odd and $v\in\mathcal S_d(\Gamma)$ an irreducible arrangement, then $v$ and $\tau v$ are deformation equivalent if and only if $v$ is flexible.
\begin{proof}
If there exists a deformation $v(\free)$ between $v$ and $\tau v$, then $\gamma(t):=\det(v,v(t))$ is continuous and satisfies $\gamma(0)=1$ as well as $\gamma(1)=-1$.
By the intermediate value theorem, there is a $t\in(0,1)$, so that, say, $\gamma(t)=1/2$. 
By \cref{res:det_properties}, the arrangement spaces of $v$ and $v(t)$ are then non-orthogonal, and by \cref{res:rigidity_characterizations}, $v$ is flexible.

The other way around, assume that $v$ is flexible. 
By \cref{res:rigidity_characterizations}, there exists a non-equivalent arrangement $\bar v$, so that $v$ and $\bar v$ have isomorphic representations. 
In fact, since $v$ and $\tau v$ have the same arrangement space, and $d$ is odd, it follows from \cref{res:det_properties} and \cref{res:representations_parity} that the representation of $\bar v$ must be \emph{positively} isomorphic to both, the representations of $v$ and $\tau v$.
With \cref{res:non_equivalent_deformation}, we then obtain a deformation $v\longrightarrow \bar v\longrightarrow \tau v$.
\end{proof}
\end{proposition}

Hence, if $v\in\mathcal S_d^+(T)$ is rigid, the open set $\mathcal S_d^+(T)$ cannot be path-connected in odd dimensions.

The following theorem is the equivalent to \cref{res:even_positively_isomorphic_iff_deformation_equivalent} in odd dimensions:

\begin{theorem}\label{res:odd_positively_isomorphic_iff_deformation_equivalent}
Let $d$ be odd, and $v,\bar v\in\mathcal S_d(\Gamma)$ irreducible and flexible arrangements. Then $v,\bar v$ are deformation equivalent if and only if they have isomorphic representations.
\begin{proof}
Again, one direction was shown in \cref{res:d_eq_implies_iso_plus}, and it remains to show that $\mathcal S_d(\Gamma)$ is path-connected.

If $v,\bar v$ have isomorphic representations, then since $d$ is odd, \cref{res:representations_parity} tells us that they also have positively isomorphic representations.
If $v$ and $\bar v$ are non-equivalent, they are then deformation equivalent by \cref{res:non_equivalent_deformation}.
If they are posi\-tive\-ly equivalent, we can already conclude deformation equivalence by the proven implication $\mathrm{eq}^+\longrightarrow\text{d-eq}$.
If, however, $v$ and $\bar v$ are equivalent, but not positively equivalent, then $\tau v$ and $\bar v$ are positively equivalent, and by $\mathrm{eq}^+\longrightarrow\text{d-eq}$ and \cref{res:odd_mirror_deformation} there exists a chain of deformations $v\longrightarrow \tau v\longrightarrow \bar v$. Hence, $v$ and $\bar v$ are deformation equivalent.
\end{proof}
\end{theorem}

The promised statement from the introduction follows as a corollary of the previous theorem, \cref{res:rigidity_characterizations} and \cref{res:non_ortho_implies_isomorphy}.

\begin{corollary}
If $d$ is odd and $v,\bar v\in\mathcal S_d(\Gamma)$ are non-equivalent, irreducible arran\-gements with non-orthogonal arrangement spaces, then $v,\bar v$ are deformation equivalent.
\end{corollary}

With the help of all of these results we can now verify the remaining implications in \cref{fig:relations}, and show that the conditions are necessary.

We already know that $\mathrm{iso}^+\longrightarrow \text{d-eq}$ holds under the assumption of even dimen\-sion. \Cref{res:odd_positively_isomorphic_iff_deformation_equivalent} now explains that the same implication is possible by just assuming flexibility (independent of the parity of the dimension).
On the other hand, if $v$ is rigid of odd dimension, then $v$ and $\tau v$ form a counterexample.

The implication $\mathrm{iso}^+\longrightarrow \mathrm{eq}^+$ (assuming even dimension \emph{and} rigidity) is obtained by taking a detour over d-eq:
$$\mathrm{iso}^+\longrightarrow \text{d-eq} \longrightarrow \mathrm{eq}^+,$$
where we need even dimension for the first implication, and rigidity for the second.
The counterexample for odd dimension is, as above, $v$ and $\tau v$.
The counterexample for flexible $v$ is obtained by taking any other non-equivalent arrangement $\bar v$ that is obtained by deforming $v$.

Finally, the implication $\mathrm{iso}\longrightarrow \mathrm{eq}$ assuming rigidity is proven as follows: if $v$ and $\bar v$ have isomorphic representations, then the representation of $v$ is positively isomorphic to either a representation of $\bar v$ or $\tau \bar v$. 
Now, either $v,\bar v$ are equivalent, or by \cref{res:non_equivalent_deformation} $v$ can be deformed into $\bar v$ or $\tau\bar v$.
But since $v$ is rigid, $v$ must already be equivalent to one (and then both) of $\bar v$ and $\tau\bar v$. 
Obvious counterexamples are obtained by taking $v$ and a non-equivalent arrangement $\bar v$ with isomorphic representation, which exists when $v$ is flexible.

Finally, we can prove an upper bound on the number of equivalence classes \wrt\ some of the relations in \cref{fig:relations}:

\begin{theorem}
For some $d\ge 1$, denote by $S_d$ the set of normalized irreducible $\Gamma$-arrangements of dimension at least $d$. It holds
\begin{myenumerate}
\item \label{it:non_isomorphic}
There are at most $n/d$ arrangements in $S_d$, so that no two have isomorphic representations.
\item \label{it:non_equivalent} 
There are at most $n/d$ rigid arrangements in $S_d$, so that no two are equivalent.
\end{myenumerate}
\begin{proof}
Let $v\^i\in S,i\in\mathcal I$ be a family of arrangements, and let $U_i\subseteq\RR^n$ be their arrangement spaces. Recall that $\dim U_i$ equals the dimension of $v\s^i$ (normalized arrangements are full-dimensional), hence $\dim U_i\ge d$. 

If $U_i\not\perp U_j$, then $v\s^i$ and $v\s^j$ have isomorphic representations by \cref{res:non_ortho_implies_isomorphy}. 
In conclusion, for \ref{it:non_isomorphic} we need to assume $U_i\perp U_j$ for all distinct $i,j\in\mathcal I$. We obtain
$$\bigoplus_{i\in\mathcal I} U_i \subseteq\RR^n \quad \implies\quad d\cdot|\mathcal I|  \le \sum_{i\in\mathcal I} \dim U_i \le n.$$
Rearranging gives the desired inequality $|\mathcal I|\le n/d$ of \ref{it:non_isomorphic}.

If now the $v\s^i$ are rigid and non-equivalent, we again obtain that their arrangement spaces are necessarily pair-wise orthogonal by \cref{res:rigidity_characterizations}, and above reasoning applies to show \ref{it:non_equivalent}.
\end{proof}
\end{theorem}

\section{Conclusion and outlook}
\label{sec:future}

We applied the concept of arrangement spaces in the context of symmetric point arrangements. 
We have shown that symmetric arrangements are in a certain one-to-one correspondence with invariant subspaces of $\RR^n$. 
We applied this to show how arrangement spaces can be used to make statements about rigidity and continuous deformations of arrangements.
In particular, we have investigated the question when a symmetric point arrangement can be continuously deformed into its mirror image without loosing symmetry in the process.

Given some $\Gamma\subseteq\Sym(N)$, we previously noted that it is in general a non-trivial task to determine the $\Gamma$-invariant subspaces of $\RR^n$. It is therefore non-trivial to obtain symmetric arrangements in this way.
However, if the point set is equipped with the additional structure of a graph, we have easy access to some of the invariant subspaces.
In fact, let $G=(N,E)$ be a graph on $N$, and $A$ its adjacency matrix.
We can now consider arrangements for which their arrangements space is an eigenspace of $A$.
These are interesting for at least two reasons:
\begin{itemize}
\item If $\Gamma\subseteq\Aut(G)$ is a subgroup of the automorphism group of $G$, then each eigenspace is already $\Gamma$-invariant. We therefore obtain an efficient tool for constructing symmetric arrangements (or better, symmetric \emph{graph realizations}), as eigenspaces are comparatively easy to compute and the computations usually also provide an appropriate ONB.
\item
The points of such an arrangement are in a specific balanced configuration. In fact, each \emph{vertex} $i\in N$ is in the span of the barycenter of its neighbors $\mathcal N_i:=\{j\in N\mid j\sim i\}$ (if the barycenter is non-zero).
\end{itemize}
%
%
More generally, such \enquote{spectral arrangements} provide a common language for ideas in graph drawing \cite{koren2005drawing}, balanced arrangements \cite{cohn2010point}, sphere packings \cite{chen2017ball} and polytope theory (specifically, in the study of \emph{eigenpolytopes}) \cite{godsil1998eigenpolytopes, padrol2010graph}
%

Already quite some work was done for such arrangements, especially in the context of distance-regular and strongly regular graphs.
Such graphs with high combinatorial regularity do not necessarily have a lot of symmetries, and hence, these classes have to be distinguished from, \eg, arc- and distance-transitive graphs.
In future publications we investigate the algebraic and geometric properties of arrangements constructed from highly symmetric graphs \cite{winter2020symmetric}, and their relation to symmetric polytopes via eigenpolytopes \cite{winter2020eigenpolytopes}.

\par\bigskip
\parindent 0pt
\textbf{Acknowledgements.} The author gratefully acknowledges the support by the funding of the European Union and the Free State of Saxony (ESF).
Furthermore, the author thanks the anonymous referee who has pointed out the close connection to the subject of frame theory and has provided several helpful references in this respect.


\bibliographystyle{abbrv}
\bibliography{literature}

\end{document}